\documentclass[11pt, letterpaper]{article}
\usepackage[utf8]{inputenc}
\usepackage{geometry}
\usepackage{csquotes}
\usepackage[spanish, english, es-nodecimaldot]{babel}

\usepackage[backend=bibtex8, citestyle=authoryear, natbib, sorting=nyt]{biblatex}
\bibliography{RandomForestsSimulation}

\usepackage{amsthm,amsmath,amsfonts,amssymb}
\usepackage{authblk}
\usepackage{graphicx}
\usepackage{float}
\usepackage{url}

\usepackage[colorlinks,citecolor=blue,urlcolor=blue,breaklinks]{hyperref}

\usepackage{algorithm}
\usepackage{algpseudocode}
\algrenewcommand{\algorithmiccomment}[1]{\textcolor{blue}{/*#1*/}}

\usepackage{nicefrac}
\usepackage{multirow}

\usepackage[english, noabbrev]{cleveref}
\usepackage{relsize}
\usepackage{bbm}
\usepackage[x11names,table]{xcolor}
\usepackage{caption}

\let\tinymatrix\smallmatrix

\patchcmd{\tinymatrix}{\vcenter}{\vtop}{}{}
\patchcmd{\tinymatrix}{\bgroup}{\bgroup\scriptsize}{}{}

\def\X{\mathbf{X}}
\def\Y{\mathbf{Y}}

\def\M{\mathbf{M}}
\def\x{\mathbf{x}}
\def\y{\mathbf{y}}

\def\imiss{\mathbf{i}}

\def\R{\mathbb{R}}          
\def\P{\mathbb{P}}          
\def\E{\mathbb{E}}          

\def\sample{\mathcal{D}_n}
\def\sampletl{\mathcal{D}_{n,\ell}}

\newcommand{\indicadora}[2]{\mathlarger{\mathbbm{1}}\,_{#1}#2}

\newcommand{\corchet}[1]{\left\lbrack #1 \right\rbrack}     
\newcommand{\parent}[1]{\left( #1 \right)}                                
\newcommand{\abs}[1]{\lvert#1\rvert}										

\DeclareMathOperator*{\argmax}{arg\,max}

\makeatletter
\def\keywords{\xdef\@thefnmark{}\@footnotetext}
\makeatother

\begin{document}

\title{Regression with Missing Data, a Comparison Study of Techniques Based on Random Forests}

\author[1,*]{Irving Gómez-Méndez}
\author[2]{Emilien Joly}
\affil[1]{Centro de Investigación en Matemáticas, AC (CIMAT)}
\affil[*]{Corresponding author: \url{irving.gomez@cimat.mx}}
\affil[2]{Centro de Investigación en Matemáticas, AC (CIMAT), \url{emilien.joly@cimat.mx}}

\date{ }
\maketitle

\begin{abstract}
In this paper we present the practical benefits of a new random forest algorithm to deal with missing values in the sample. The purpose of this work is to compare the different solutions to deal with missing values with random forests and describe our new algorithm performance as well as its algorithmic complexity. A variety of missing value mechanisms (such as MCAR, MAR, MNAR) are considered and simulated. We study the quadratic errors and the bias of our algorithm and compare it to the most popular missing values random forests algorithms in the literature. In particular, we compare those techniques for both a regression and prediction purpose. This work follows a first paper \citet{gomezmendez2020consistency} on the consistency of this new algorithm.
\keywords{\emph{Mathematics Subject Classification 2020.} Primary 62G08; Secondary 60G25}
\keywords{\emph{Key words and phrases.} Missing values, Random forests, Non-parametric regression, Prediction with missing values}
\end{abstract}

\section{General introduction}

\subsection{Random forests with missing values: a contemporary challenge}
Random forests and recursive trees are widely used in applied statistics and computer science. The popularity of recursive trees relies on several factors: their easy interpretability, the fact that they can be used for both regression and classification tasks, the small number of hyper-parameters to be tuned and finally, their non-parametric nature that allows their use to infer arbitrarily complex relations between the input and the output space. A random forest combines several randomized trees, improving the prediction accuracy at a cost of a slight lost in interpretation. This technique is easily parallelizable which has made it one of the most popular tools for handling high dimensional data sets. It has been successfully involved in various practical problems, including chemioinformatics, ecology, 3D object recognition, bioinformatics and econometrics. \citet{biau2016random} present a detailed list of applications as well as a review on random forests. In the present work we have focused on the ability of random forests to deal with missing values.

\paragraph{Three groups of algorithms.}
Approaches to handle missing values regarding the use of random forests could be classified in three groups. In the \textit{first group}, we classify the algorithms that impute the missing values without using ad-hoc techniques, like $k$ nearest neighbor imputation ($knn$-imputation) \citep{troyanskaya2001missing} or multiple imputation chained equations (MICE) \citep{van2006fully}, for example. Once performed this imputation step, the regular random forest algorithm is implemented on the completed data set. In this paper, we will not pay much attention to this first group of algorithms as the imputation step is unrelated with the later use of the random forest algorithm which is finally, only treated as a black box.
The \textit{second group} is composed of the algorithms that use the random forest structure (directly or through some extra method like proximity matrices) to impute the missing values. They typically use iterated updates of the imputations to refine them. Logically, we decided to store the methods proposed by \citet{Breiman2003,ishioka2013imputation,stekhoven2011missforest} in this group. 
The \textit{third group} consists in built-in methodologies that include the incomplete observations directly in the construction of the recursive trees (without imputation steps) to compute the random forest's estimation. Are included in this group: surrogate splits \citep{breiman1984CART}, missing incorporated in attributes (MIA) \citep{twala2008good} or the algorithm proposed by \citet{gomezmendez2020consistency}, which is the main focus of this paper. 

\subsection{Literature on missing values treatment with random forests}

Handling missing values through different algorithms has recently received much attention with several simulation studies comparing the performance of distinct methods \citep{feelders1999handling,farhangfar2008impact,rieger2010random,hapfelmeier2012recursive,josse2019consistency}. In this paper, we propose a comparison study of techniques of the last two groups since they appear to be the most used and precise in practice. Of special focus for this present paper is the study of the practical performance of the new approach presented in \citet{gomezmendez2020consistency}, which has the advantage to be consistent under an MCAR mechanism and a generalized additive model.

The idea of handling missing values with random forests algorithms is not new whether it be for a ``filling the gaps'' task, a learning or a prediction task with corrupted data. In this section, we describe the current state of the art in the study of this problem.
The study developed by \citet{feelders1999handling} is one of the first to present a result on missing values using recursive trees. This work compares the error rate between surrogate splits in a single decision tree, and the imputation procedure presented by \citet{schafer1997analysis,schafer1998multiple}. Two data sets are considered in the study, given by the so-called \textit{Waveform} data set (presented by \citet{breiman1984CART}) and the \textit{Prima Indian Datasets} (available at the UCI machine learning repository). The percentage of missing values varies between 10\% and 45\% for the first data set, where the missing values are introduced accordingly to an MCAR mechanism. On the other hand, 10\% of the observations present missing values in the second data set. \citet{feelders1999handling} concludes that the imputation procedure yields significantly less missclassification error rate.

The work of \citet{farhangfar2008impact} presents a comparison of classification methods like support vector machines, $knn$-imputation and the decision tree method C4.5 \citep{quinlan1993} applied to missing data, and imputation algorithms, including single imputation and MICE. In total 15 data sets are considered, inducing missing values with up to 50\% of the observations per variable being missing. The authors conclude that the application of MICE leads to better results in most of the instances.

\citet{rieger2010random} compare surrogate splits introduced in the conditional inference forests \citep{hothorn2006unbiased} with $knn$-imputation, with no clear advantage for either of the methods. This study considers classification and regression problems with three different correlation structures and seven schemes for missing values. However, the percentage of missing values is kept constant.

\citet{hapfelmeier2012recursive} present a comparison study using trees built with the CART criterion, conditional inference trees and their corresponding random forests, focusing on surrogate splits to handle missing values and the use of MICE to impute missing values. The authors consider 12 real life data sets, half of them for regression and the other half for classification, 8 of the data sets already present missing values while in the rest 4 data sets missing values are induced. Arguing that \citet{rieger2010random} found similar results for MCAR and MAR mechanisms, the missing values are solely introduced according to an MCAR missing-data mechanism, considering missing rates between 0\% (benchmark) to 40\%. Their results do not show a clear improvement by using multiple imputation, with MICE even producing inferior results when missing values are limited in number and are not arbitrary spread across the data. The results also show a similar result between trees constructed with the CART criterion and conditional inference trees.

The simulation developed by \citet{josse2019consistency} studies the performance of several methods to handle missing values in regression tasks. 
In this study, three different regression functions are considered, one of them being linear, one quadratic and the third being the so-called ``friedman1''  \citep{friedman1991multivariate}. 
They consider the MCAR mechanism and censoring, inducing up to 20\% of observations with missing values in the first variable. 
The authors compare conditional inference trees, trees and random forests based on the CART criterion and XGBoost \citep{chen2016xgboost}. 
The algorithms to handle missing values include MIA, surrogate splits for both trees built with the CART criterion and conditional inference trees, block propagation (only implemented in XGBoost), surrogate splits, mean-imputation and EM imputation \citep{dempster1977maximum}. 
The authors clearly favors the usage of MIA for tree-based methods, while block propagation could also be a good method.

\citet{Breiman2003} proposes an algorithm based on random forests to impute the missing values making use of the observed values and the proximity matrix.
\citet{ishioka2013imputation} presents an improvement which instead of considering only the observed values in the imputation procedure, it considers nearest neighbors for continuous variables and all the observations for categorical variables. \citet{ishioka2013imputation} compares these two approaches and $knn$-imputation introducing missing values completely at random in the \textit{Spam} data set and considering missing data rates from 5\% to 60\%, concluding that the proposed approach outperforms $knn$-imputation and the previous method proposed by \citet{Breiman2003}. However, the same author comments that other missing-data mechanisms should be considered.

\citet{stekhoven2011missforest} introduce the missForest algorithm which imputes the missing values iteratively considering it as a regression problem in which the imputation of the current variable is done using all the other variables. In the simulation study presented, they consider classification and regression problems in 7 different data sets where 10\%, 20\% or 30\% of the values are removed completely at random, concluding that missForest outperforms $knn$-imputation and MICE. Furthermore, the authors ensure that the full potential of missForest is deployed when data includes interactions or non-linear relations between variables of unequal scales and different types.

In \Cref{sec:results}, we compare the performance of most of these techniques with our new proposal taken from \citet{gomezmendez2020consistency} when we let the percentage of missing values to vary from 0\% to 95\%. 

The rest of the paper is organized as follows. In \Cref{sec:concepts}, we give the general background on random forests algorithms. We formally introduce the CART criterion and the missing-data mechanisms. In \Cref{sec:newapproach}, we describe our proposal to build recursive trees including the missing values and we give upper bounds on its computational complexity. In \Cref{sec:related_methods}, we describe the other methods that we use as benchmark in the simulation study explained in detail in \Cref{sec:simulations}. In \Cref{sec:results}, we presents and describe the results. Finally \Cref{sec:conclusions} presents the conclusions.

\section{Settling the concepts}
\label{sec:concepts}

\subsection{Random forests built upon the CART criterion}

Throughout this article, we assume to have access to a training data set $\sample=(\X_i,Y_i)_{i=1,\dots,n}$ where the response variables $Y_i$ are real-valued and the input variables $\X_i$ belong to some space $\mathcal{X}\subseteq \R^{p}$. The objective is to use the data $\sample$ to construct a learning model, also called learner, predictor or estimator, $m_n:\mathcal{X}\to\R$ that estimates the regression function $m(\x)=\E[Y|\X=\x]$.

The random forest is made of a set of regression trees that are later aggregated all together with a simple mean idea.
Each branch of the tree will be random (in its construction process) and represents a partition of the input space in smaller regions.
Moving along the path of the tree corresponds to a choice of one of the possible regions.
To construct this partition of the input space, the trees are built in a recursive way (hence the name of recursive trees). 
The root of the tree corresponds to the whole input space $\mathcal{X}$. 
Then, recursively, a region is chosen and is split into two smaller regions. This process is continued until some stopping rule is applied. 
At each step of the tree construction, the partition performed over a cell (or equivalently its corresponding node) is determined by maximizing some split-criterion. 
The present work focuses in the so-called CART split criterion.
We first introduce some important notations.

\begin{itemize}
\item $A$ denotes a general node (or cell).

\item $N(A)$ holds for the number of points in $A$.

\item The notation $d= (h,z)$ denotes a cut in $A$, where\medskip

$h$ is a direction, $h\in\{1,\ldots,p\}$, and

$z$ is the position of the cut in the $h$th direction, between the limits of $A$.

\item $\mathcal{C}_A$ is the set of all possible cuts in node $A$.

\item A cell $A$ is split into two cells denoted $A_L=\{\x\in A\,:\,\x^{(h)}< z\}$ and $A_R=\{\x\in A\,:\,\x^{(h)}\geq z\}$.

\item $\bar{Y}_A$ (resp. $\bar{Y}_{A_L}$, $\bar{Y}_{A_R}$) is the empirical mean of the response variable $Y_i$ for the indexes such that $\X_i$ belongs to the cell $A$ (resp. $A_L$, $A_R$).
\end{itemize}
Then, the CART split criterion for a generic cell $A$ is defined as
\begin{align}
L_n(A,d)=&\frac{1}{N(A)}\sum_{i=1}^n\parent{Y_i-\bar{Y}_A}^2\indicadora{\X_i\in A}{}\nonumber\\
&-\frac{1}{N(A)}\sum_{i=1}^n\left(Y_i-\bar{Y}_{A_L}\indicadora{\X_i^{(h)}<z}{}-\bar{Y}_{A_R}\indicadora{\X_i^{(h)}{}\geq z}\right)^2\indicadora{\X_i\in A}{}\label{eq:CART_definition}
\end{align}
with the convention $0/0=0$.

As mentioned above, a random forest is a predictor consisting of $M(>1)$ randomized trees.
The randomization is introduced in two different parts of the tree construction. 
Prior to the construction of each tree, $a_n$ observations are extracted at random with (or without) replacement from the learning data set $\sample$.
Only these $a_n$ observations are taken into account in the tree construction. 
Then, at each cell a split (or cut) is performed by maximizing the split criterion over a number $\texttt{mtry}$ of directions $h$, chosen uniformly at random. The tree construction is stopped when each final node contains less or equal than \texttt{nodesize} points. Hence, the parameters of this algorithm are:

\begin{itemize}
\item $M>1$, which is the number of trees in the forest.

\item $a_n\in\{1,\ldots,n\}$, which is the number of observations in each tree.

\item $\texttt{mtry}\in\{1,\ldots,p\}$, which is the number of directions (features) chosen, candidates to be split. We denote by $\mathcal{M}_{try}$ the features selected in each step.

\item $\texttt{nodesize}\in\{1,\ldots,a_n\}$, which is the maximum number of observations for a node to be a final cell.
\end{itemize}

The randomization introduced in the trees (independent from the original source of randomness in the sample \(\sample\)) is represented in a symbolic random variable \(\Theta\). To each tree -- randomized with the random variable \(\Theta_k\) -- there is associated a predicted value at a query point $\x$, denoted as $m_n(\x;\Theta_k)$. The different trees are constructed by the same procedure but with independent randomization, so the random variables $\Theta_1,\ldots,\Theta_M$ are i.i.d. with common law $\Theta$. In our choice of the construction rules, \(\Theta\) consists in the observations selected for the tree and the candidate variables to split at each step. Finally, the $k$th tree's estimation at point $\x$ is defined as
\[
m_n(\mathbf{x};\Theta_k)=\sum_{i\in \mathcal{I}_{n,\Theta_k}}\frac{Y_i\indicadora{\X_i\in A_n(\x;\Theta_k)}{}}{N(A_n(\x;\Theta_k))}
\]
where $\mathcal{I}_{n,\Theta_k}$ is the set of the \(a_n\) observations selected prior to the construction of the $k$th tree, $A_n(\x;\Theta_k)$ is the unique final cell that contains $\x$, and $N(A_n(\x;\Theta_k))$ is the number of observations which belong to the cell $A_n(\x;\Theta_k)$. The average of the trees forms the random forest's estimation given by

\[ m_{M,n}(\x;\Theta_1,\ldots,\Theta_M)=\frac{1}{M}\sum_{k=1}^M m_n(\x;\Theta_k).\]

It is known from the work of \citet{breiman2001random} that the random forest does not overfit when $M$ tends to infinity. This makes the parameter \(M\) only restricted by computational power.

\subsection{Missing-data mechanisms}

The concept of missing-data mechanism (introduced by \citet{rubin1976inference}) establishes the relationship between missingness and data. Before introducing the missing-data mechanisms,  let us define a new variable, called the missing-data indicator
\[\M^{(h)}=\left\{\begin{array}{ll}
1 & \text{if } \X^{(h)}\text{ is missing}\\
0 & \text{otherwise}
\end{array} \right.,\quad 1\leq h\leq p.\]
We assume throughout this work that the response $Y$ has no missing values which makes unnecessary to define an indicator of missing variable for $Y$. Then, the mechanisms are fully characterized by the information of the conditional distribution of $\M^{(h)}$ given $(\X,Y)$. There are three possible missing-data mechanisms.

\paragraph{Missing Completely at Random (MCAR).} We say that the variable $\X^{(h)}$ is MCAR if $\M^{(h)}$ is independent from $(\X,Y)$. In other words, under the MCAR assumption, a coordinate $\X^{(h)}$ has some probability to be missing in the sample and this probability does not depend on the value of \(\X\) nor the response variable \(Y\).

\paragraph{Missing at Random (MAR).} Let us define the set $h_o=\{h:\P\corchet{\M^{(h)}=0}=1\}$, thus $\X^{(h_o)}$ is the vector formed by the variables that are always observed.
The variable $\X^{(h)}$ is MAR if  \[\P\corchet{\M^{(h)}=1|(\X,Y)}=\P\corchet{\M^{(h)}=1|(\X^{(h_o)},Y)}\] That is, the probability of $\X^{(h)}$ being missing only depends on observed data.

\paragraph{Missing Not at Random (MNAR).} If the probability of missingness depends on unobserved values we say that $\X^{(h)}$ is MNAR.

\section{A new approach}
\label{sec:newapproach}
It is of special interest to study the performance of the algorithm presented by \citet{gomezmendez2020consistency} who have proven its consistency for an MCAR mechanism and a generalized additive model, being one of the first results on the consistency of random forests with missing values. This algorithm adapts the original CART criterion to manage missing data directly in the construction of the regression trees. We now recall this proposal.

\subsection{Description of the algorithm}

Assume for now that we have a partition of the input space. When there are missing values in the data set, there is uncertainty on the region to which each observation belongs. Thus, the original CART criterion becomes intractable since the quantities $N(A)$, $N(A_L)$, $N(A_R)$, $\bar{Y}_A$, $\bar{Y}_{A_L}$, $\bar{Y}_{A_R}$, $\indicadora{\X_i\in A}{}$, $\indicadora{\X_i^{(h)}<z}{}$ and $\indicadora{\X_i^{(h)}\geq z}{}$ can not be computed. The proposed approach keeps the form of the CART criterion and makes use of adapted imputations for the intractable parts which allows the computation of a modified version of the CART criterion. Unlike most of the imputation techniques, the imputation step is not performed independently of the evaluation of the CART criterion but is integrated to its later optimization. While a cut is selected by maximizing the original CART criterion, now a couple (cut, imputation) is chosen at each split in the creation of the random tree. The idea is that, for a cut, the observations with missing values are assigned to the child node that maximizes this modified CART criterion. At the end, the missing observations will belong to a final node of the tree, which can give an ``imputation'' of the missing values as a region of the input space, which in turn would be translated into a ``cloud'' of possible regions for the missing values when the random forest is considered. For seek of clarity the proposal has been divided in two parts, \Cref{alg:random_forests_with_missingness} describes the steps for the construction of the random forest with missing entries and \Cref{alg:CART_with_missingness} which establishes the steps to find the best cut and assignation of missing data. At the beginning of the construction of each tree it is necessary to initialize the list of current not-final cells $\mathcal{P}=\{\mathcal{X}\}$ as well as the induced partitions of the input variables and the target variable $\mathcal{X}_{\mathcal{P}}=\{(\X_1,\ldots,\X_n)\}$ and $\mathcal{Y}_{\mathcal{P}}=\{(Y_1,\ldots,Y_n)\}$, respectively. Along with the initialization of these sets, the final partition of the input space $\mathcal{P}_f=\{\}$ and the corresponding partitions of the input variables $\mathcal{X}_f=\{\}$ and the target variable $\mathcal{Y}_f=\{\}$ are initialized. When a cell $A$ satisfies the criteria to be a final cell, it is removed from $\mathcal{P}$ and added to $\mathcal{P}_f$, the input variables $\mathcal{X}_A$ and the target variables $\mathcal{Y}_A$ belonging to $A$ are also removed from $\mathcal{X}_{\mathcal{P}}$ and $\mathcal{Y}_{\mathcal{P}}$, and added to $\mathcal{X}_{f}$ and $\mathcal{Y}_{f}$, respectively. In the sequel, we denote as $N_{obs}^{(h)}(A)$ the number of points belonging to $A$, whose value in the direction $h$ has been observed and $N_{miss}^{(h)}(A)$ the number of observations assigned to cell $A$ whose value in the direction $h$ is missing.

\begin{algorithm}[H]
\hspace*{\algorithmicindent} \textbf{Input:} Training sample $\sample$, number of trees $M>1$, \texttt{mtry}$\in\{1,\ldots,p\}$, $a_n\in\{1,\ldots,n\}$, \texttt{nodesize}$\in\{1,\ldots,a_n\}$.\\
\hspace*{\algorithmicindent} \textbf{Output:} Random forest $m_{M,n}$.
\begin{algorithmic}[1]
\For{$i=1,\ldots,M$}
	\State Select $a_n$ points uniformly in $\sample$.
	\State Initialize $\mathcal{P}$, $\mathcal{X}_{\mathcal{P}}$, $\mathcal{Y}_{\mathcal{P}}$, $\mathcal{P}_f$, $\mathcal{X}_f$ and $\mathcal{Y}_f$.
	\While{$\mathcal{P}\neq \varnothing$}
		\State Let $A$, $\mathcal{X}_A$, $\mathcal{Y}_A$ be the first elements of $\mathcal{P}$, $\mathcal{X}_{\mathcal{P}}$ and $\mathcal{Y}_{\mathcal{P}}$, resp.
		\State Set $N(A)$ the number of points which belong or were assigned to $A$.
		\If{$N(A)\leq$\texttt{nodesize}}
			\State Remove $A$, $\mathcal{X}_A$ and $\mathcal{Y}_A$ from $\mathcal{P}$, $\mathcal{X}_{\mathcal{P}}$ and $\mathcal{Y}_{\mathcal{P}}$, and add them to $\mathcal{P}_f$, $\mathcal{X}_f$ and $\mathcal{Y}_f$.
		\Else
			\For{$j=1,\ldots,p$}
				\State Compute $N_{obs}^{(j)}(A)$.
			\EndFor
			\State Let $h_{obs}$ be the features $h$ such that $N_{obs}^{(h)}(A)>1$.
			\If{$h_{obs}=\varnothing$}
				\State Remove $A$, $\mathcal{X}_A$ and $\mathcal{Y}_A$ from $\mathcal{P}$, $\mathcal{X}_{\mathcal{P}}$ and $\mathcal{Y}_{\mathcal{P}}$, and add them to $\mathcal{P}_f$, $\mathcal{X}_f$ and $\mathcal{Y}_f$.
			\Else
				\State Set $m_{obs}=|h_{obs}|$.
				\If{$m_{obs}\leq$\texttt{mtry}}
					\State Set $\mathcal{M}_{try}=h_{obs}$.
				\Else
					\State Select uniformly, without replacement, a subset $\mathcal{M}_{try}\subset h_{obs}$ of cardinality \texttt{mtry}.
				\EndIf
				\State Apply \Cref{alg:CART_with_missingness} on cell $A$ along the features in $\mathcal{M}_{try}$.
				\State Call $A_L$ and $A_R$ the two resulting cells.
				\State Remove $A$, $\mathcal{X}_A$ and $\mathcal{Y}_A$ from $\mathcal{P}$, $\mathcal{X}_{\mathcal{P}}$ and $\mathcal{Y}_{\mathcal{P}}$.
				\State Add $A_L$, $A_R$, $\mathcal{X}_{A_L}$, $\mathcal{X}_{A_R}$, $\mathcal{Y}_{A_L}$ and $\mathcal{Y}_{A_R}$ to $\mathcal{P}$, $\mathcal{X}_{\mathcal{P}}$ and $\mathcal{Y}_{\mathcal{P}}$.
			\EndIf
		\EndIf
	\EndWhile
\EndFor
\end{algorithmic}
\caption{Random forest with assignation of missing entries.} \label{alg:random_forests_with_missingness}
\end{algorithm}

\begin{algorithm}[H]
\hspace*{\algorithmicindent} \textbf{Input:} Cell $A$, $(\mathcal{X}_A,\mathcal{Y}_A)$, $\mathcal{M}_{try}$.\\
\hspace*{\algorithmicindent} \textbf{Output:} Best cut and assignation $(\hat{z},\hat{w})$.
\begin{algorithmic}[1]
	\State Set $max_{CART}\leftarrow 0$
	\For{$h\in\mathcal{M}_{try}$}
		\State Compute the midpoint of two consecutive values of $\X^{(h)}$ between those points $\X\in A$, let be $Z_A^{(h)}$ the set of these midpoints.
		\State Let be $W_A^{(h)}$ the set with all the possible assignation for the missing values in $h$.
		\For{$z\in Z_A^{(h)}$}
			\For{$w\in W_A^{(h)}$}
				\State Let $c_A$ be the CART-criterion computed with the cut $(h,z)$ and the assignation $w$.
				\If{$c_A > max_{CART}$}
					\State $max_{CART}\leftarrow c_A$.
					\State $(\hat{z},\hat{w})\leftarrow (z,w)$
				\EndIf
			\EndFor
		\EndFor
	\EndFor
\end{algorithmic}
\caption{Best cut and assignation.} \label{alg:CART_with_missingness}
\end{algorithm}

\subsection{Complexity of the Algorithm}
\label{sec:algo_complexity}

At first sight, it looks that the number of possible assignations defined by $\mathcal{W}_A^{(h)}$ in Algorithm \ref{alg:CART_with_missingness} are all the combinations for the assignations of the $N^{(h)}_{miss}(A)$ missing observations. This would lead to exponential complexity in the number of calculations of the CART criterion and then to an intractable algorithm. However, the number of possible candidate assignations to maximize the CART criterion is lower than this quantity. To see this, fix a cell $A$ and a cut $(h,z)$ in $A$, to keep a simple notation let $\bar{Y}_{L,obs}$ (resp. $\bar{Y}_{R,obs}$) be the mean of the response variable for the points belonging to the left (right) node and observed in the direction $h$. Suppose without lost of generality that $\bar{Y}_{L,obs}\leq \bar{Y}_{R,obs}$ and denote by $\imiss_{miss}=\{1,\ldots,N\}$ the set of indexes of the observations assigned to the cell $A$ whose direction $h$ is missing, without lost of generality assume that $Y_{1}\leq\cdots\leq Y_{N}$. Because $\bar{Y}_{L,obs}\leq\bar{Y}_{R,obs}$ and maximizing the CART criterion implies to make as different as possible the average of the target on the left node from the average of the target on the right node, then  observations $i\in\imiss_{miss}$ with the lowest values $Y_{i}$ should be assigned to the left node and the observations $i\in\imiss_{miss}$ with the largest values $Y_{i}$ should be assigned to the right node. Therefore, there exists $w\in\{1,\ldots,N\}$ such that assigning
$Y_{1},\ldots,Y_{w-1}\in A_L$ and $Y_{w},\ldots,Y_{N}\in A_R$ maximizes the CART criterion. Hence, the set with all possible assignations $W_A^{(h)}$ has a cardinality of $N_{miss}^{(h)}(A)+1$ and there is only a linear number of assignations to be considered. Now, we can calculate the complexity of the algorithm.


Let be $O_{CART}$ the number of operations needed to calculate the CART criterion for a given cut $(h,z)$ and assignation $w$ of missing values. It makes sense to consider the complexity of our algorithm with respect to $O_{CART}$ since every other random forest algorithm that we compare to also make a certain number of evaluation of the CART criterion. Consider the \Cref{alg:CART_with_missingness} and note that $\abs{Z_A^{(h)}}=N_{obs}^{(h)}(A)-1$ and  $\abs{W_A^{(h)}}=N_{miss}^{(h)}(A)+1$. Thus, the number of necessary operations to the get the best cut and assignation is \[\sum_{h\in\mathcal{M}_{try}}\abs{Z_A^{(h)}}\abs{W_A^{(h)}}O_{CART}\] On the other hand, denote by $\mathcal{P}_{nf}$ the set of non-final nodes of a tree, it is clear that $|\mathcal{P}_{nf}|\leq n-1$,
where the equality holds when each final cell contains just one observation. Now, let be $O_{tree}$ the number of operations needed to build a regression tree with our approach, then

\begin{align*}
O_{tree}&=\sum_{A\in\mathcal{P}_{nf}}\sum_{h\in\mathcal{M}_{try}}\abs{Z_A^{(h)}}\abs{W_A^{(h)}}O_{CART}\\
&=\sum_{A\in\mathcal{P}_{nf}}\sum_{h\in\mathcal{M}_{try}}\parent{N_{obs}^{(h)}(A)-1}\parent{N_{miss}^{(h)}(A)+1}O_{CART}\\
&\leq \sum_{A\in\mathcal{P}_{nf}}\sum_{h\in\mathcal{M}_{try}}\parent{N(A)-1}\parent{N(A)+1}O_{CART}\\
&\leq \texttt{mtry}\times n^3\times O_{CART}\\
\end{align*}

\subsection{On simplifications of the algorithm}
\label{sec:simplifications}

In this section we discuss a remark that leads to further simplification of our algorithm in order to reach an algorithmic complexity of the same order than MIA. Recall from \Cref{sec:algo_complexity} that the algorithmic complexity of our proposal is of the order \(O(n^3)\) times the calculation required to compute one single CART criterion and that MIA algorithm (in \(O(n^2)\)) is quicker by a linear factor.

\begin{figure}[H]
\begin{center}
\includegraphics[width=\textwidth]{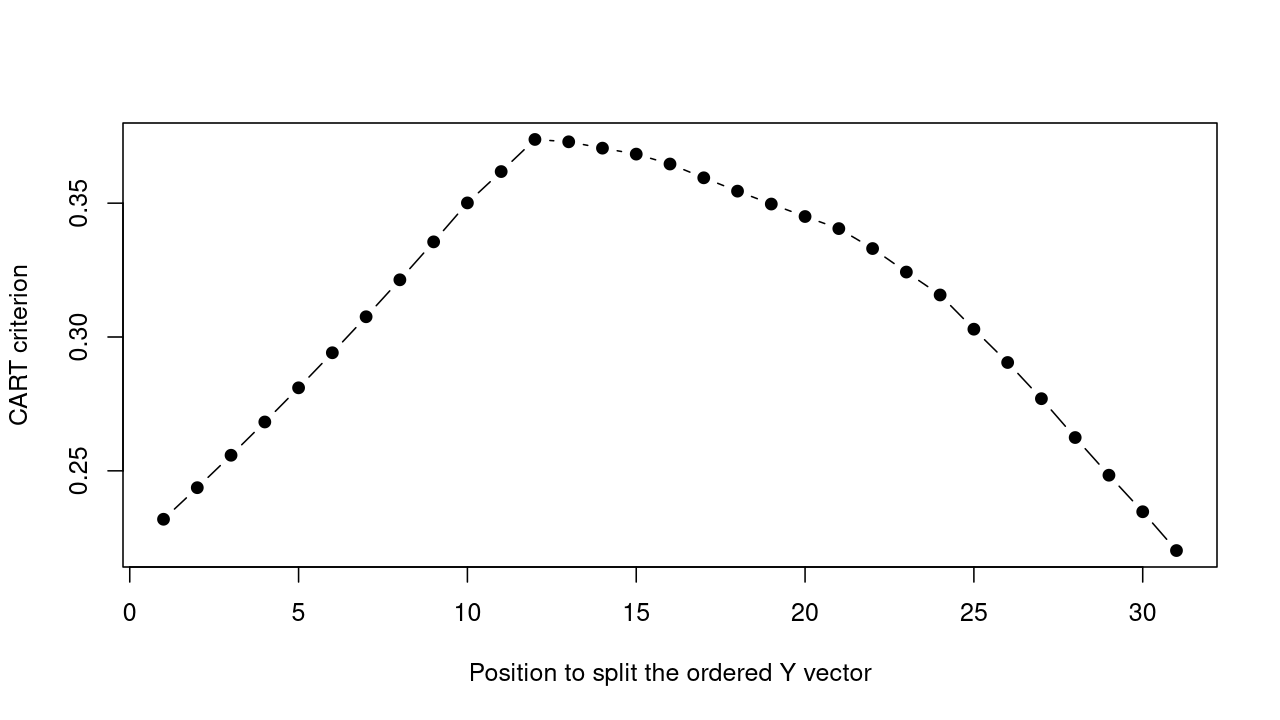}
\caption{CART criterion as a function of the position where we split $Y_{1},\ldots,Y_{N}$.}
\label{fig:CART_split_Y}
\end{center}
\end{figure}

From the CART criterion written in \Cref{eq:CART_definition}, we see that the criterion is convex with respect to the variable that counts the number of observations assigned into \(A_L\), see \Cref{fig:CART_split_Y} for a graphical representation of this fact. This simple remark leads us to opt for a dichotomy strategy for the optimal assignation of the missing variables.
Indeed, at first step, we assign half of the missing variables to the left (in \(A_L\)) and half to the right (in \(A_R\)). The corresponding assignation is then given by a \(k = \lfloor N/2 \rfloor\) (called pivot for the dichotomy) such that \(Y_1\le \dots \le Y_k \in A_L\) and \(Y_{k+1}\le \dots \le Y_N \in A_R\), assuming once again that $\bar{Y}_{L,obs}\leq \bar{Y}_{R,obs}$ without lost of generality. Such configuration gives a CART criterion resumed (abusively) in the notation \(CART(k)\). We search for the optimal assignation through calculations of the local gradient given by
\begin{equation*}
  \nabla CART(k)=CART(k+1)-CART(k).
\end{equation*}
If the value is positive (resp. negative), then the optimal $k$ is bigger or equal (resp. smaller or equal ) to $\lfloor N/2\rfloor$. Say (for example), that $\nabla CART(k)>0$. Then, we place the new pivot at the point $k=\lfloor 3N/4\rfloor$ and compute $\nabla CART(k)$. Once again, the sign of the gradient tells us in which sub interval of $[\lfloor N/2\rfloor,N]$ the optimal assignation is.
We, then repeat this simple searching procedure until ending with an interval containing only two points. The biggest $CART$ value of the two gives the optimal assignation of our variables. This dichotomy procedure allows to consider only $\log (N)$ computations of the local gradient and a simple comparison at the end. This allows us to replace the complexity of the algorithm by
\[
O_{tree}\leq \texttt{mtry}\times 2n^2\log(n)\times O_{CART}
\]
which is comparable to the MIA algorithm complexity up to a logarithmic factor. This optimization of the procedure is particularly important when the algorithm deals with a very corrupted data set where the missing entrees could represent a significant proportion of the all data.

\section{Experimental details}
\label{sec:related_methods}

\subsection{Benchmark methods}
Many methods proposed in the literature to handle missing data using random forests operate through imputation in a recursive way. They start by using the original training data set $\sample$ to fill the blank spaces in a roughly way. For example, using the median of the observed values in the specific direction. We denote this new data set as $\mathcal{D}_{n,1}$.

The imputed data set $\mathcal{D}_{n,1}$ is used to build a random forest. Then, some structures of the forest are exploited, like the so-called proximity matrix, improving the imputation and resulting in a new data set $\mathcal{D}_{n,2}$. The procedure continues by iterations until some stopping rule is applied. These stopping rules trigger, for example, when the update in imputed variables becomes negligible or when a fix number of iterations is achieved. More formally, let us define
\[\X_{i,\ell}^{(h)}=\left\{\begin{array}{ll}
\X_i^{(h)} & \text{if } \M_i^{(h)}=0\\
\widehat{\X}_{i,\ell}^{(h)} & \text{if } \M_i^{(h)}=1
\end{array} \right.\]
where $\widehat{\X}_{i,\ell}^{(h)}$ is the imputation of $\X_{i}^{(h)}$ at time $\ell\geq 1$, and let  $\X_{i,\ell}=\parent{\X_{i,\ell}^{(1)},\ldots,\X_{i,\ell}^{(p)}}$.
To properly introduce the methods considered in the simulation study, we need to define the connectivity between two points in a tree and the proximity matrix of the forest.
Let $K_{\Theta,n}(\X,\X')$ be 1 if and only if $\X$ and $\X'$ belong to the same final cell (we say that $\X$ and $\X'$ are connected in the tree $m_n(\cdot;\Theta)$) in the tree designed with $\sample$ and the parameter $\Theta$. Otherwise, $K_{\Theta,n}(\X,\X')=0$. Finally, the proximity between  between $\X$ and $\X'$ in the random forest $m_{M,n}(\cdot;\Theta_1,\ldots,\Theta_M)$, is defined as
\[K_{M,n}(\X,\X')=\frac{1}{M}\sum_{k=1}^M K_{\Theta_k,n}(\X,\X').\]
Analogously, we define the proximity $K_{M,\ell}(i,j)$ between $\X_i$ and $\X_j$ at time $\ell$ (i.e. using the data set $\sampletl$). 
We also define $\mathbf{i}_{miss}^{(h)}\subseteq\{1,\ldots,n\}$ as the indexes where $\X^{(h)}$ is missing, and $\mathbf{i}_{obs}^{(h)}=\{1,\ldots,n\}\setminus\mathbf{i}_{miss}^{(h)}$ as the indexes where $\X^{(h)}$ is observed.

We consider three different approaches which impute missing values through random forests. These methods correspond to the algorithms presented respectively in \citet{Breiman2003,ishioka2013imputation,stekhoven2011missforest}. We refer to those as Breiman's approach, Ishioka's approach and missForest approach, respectively. These algorithms impute the missing values through iterative improvements and hence belong to the second category. We also consider MIA, which is an algorithm that handle the missing values directly in the construction of the trees, assigning all the missing values to the same cell. As simple baselines we consider median-imputation of the missing values and listwise deletion (i.e. removing the observations with missing values) before the construction of the random forest.

\subsubsection{With imputation}
In this section, we describe the algorithms that belong to the second group (with iterative imputation) that are included in the upcoming simulations.
\paragraph{Breiman's Approach.}
If $\X^{(h)}$ is a continuous variable, $\widehat{\X}_{j,\ell+1}^{(h)}$ is the weighted mean of the observed values in $\X^{(h)}$, where the weights are defined by the proximity matrix of the previous random forest, that is
\[\widehat{\X}_{j,\ell+1}^{(h)}=\frac{\sum_{i\in\mathbf{i}_{obs}^{(h)}}K_{M,\ell}(i,j)\X_i^{(h)}}{\sum_{i\in\mathbf{i}_{obs}^{(h)}}K_{M,\ell}(i,j)},\quad \begin{array}{l}\ell\geq 1\\
j\in\mathbf{i}_{miss}^{(h)}
\end{array}\]
On the other hand, if $\X^{(h)}$ is a categorical variable, $\widehat{\X}_{j,\ell+1}^{(h)}$ is given by
\[\widehat{\X}_{j,\ell+1}^{(h)}=\argmax_{\x\in\mathcal{X}^{(h)}}\sum_{i\in \mathbf{i}_{obs}^{(h)}}K_{M,\ell}(i,j)\indicadora{\X_i^{(h)}=\x}{},\quad \begin{array}{l}\ell\geq 1\\
j\in\mathbf{i}_{miss}^{(h)}
\end{array}\]
That is, $\widehat{\X}_{j,\ell+1}^{(h)}$ is the class that maximizes the sum of the proximity considering the observed values in the class.

\paragraph{Ishioka's Approach.}
If $\X^{(h)}$ is a continuous variable, $\widehat{\X}_{j,\ell+1}^{(h)}$ is the weighted mean of the $k$ nearest neighbors, according to the proximity matrix, over all the values, both imputed and observed. The $k$ closest values are chosen to make more robust the method and avoid values which are outliers.
\[\widehat{\X}_{j,\ell+1}^{(h)}=\frac{\sum_{\begin{smallmatrix} i\in\text{neigh}_k \\ i\neq j
\end{smallmatrix}}K_{M,\ell}(i,j)\widehat{\X}_{i,\ell}^{(h)}}{\sum_{\begin{smallmatrix} i\in\text{neigh}_k \\ i\neq j
\end{smallmatrix}}K_{M,\ell}(i,j)},\quad \begin{array}{l}\ell\geq 1\\
j\in\mathbf{i}_{miss}^{(h)}
\end{array}\]
For categorical variables, it is not necessary to see only the $k$ closest values because the outliers of $\X$ will have few attention. Meanwhile the proximity with missing values should have more attention, especially when the missing rate is high. Hence, if $\X^{(h)}$ is a categorical variable, $\widehat{\X}_{j,\ell+1}^{(h)}$ is given by

\[\widehat{\X}_{j,\ell+1}^{(h)}=\argmax_{\x\in\mathcal{X}^{(h)}}\sum_{i\neq j}K_{M,\ell}(i,j)\indicadora{\widehat{\X}_{i,\ell}^{(h)}=\x}{},\quad \begin{array}{l}\ell\geq 1\\
j\in\mathbf{i}_{miss}^{(h)}
\end{array}\]

\paragraph{MissForest.}
This algorithm treats the imputation as a regression problem by itself, where the target variable is the input variable with missing values. MissForest predicts the missing values using a random forest trained on the observed parts of the data set. More formally for an arbitrary variable $\X^{(h)}$ we can separate the data set into four parts:

\begin{itemize}
\item the observed parts of the variable $\X^{(h)}$, denoted as $\y_{obs}^{(h)}$;
\item the missing values of the variable $\X^{(h)}$, denoted as $\y_{miss}^{(h)}$;
\item the variables other than $\X^{(h)}$ and the observations whose indexes belong to $\mathbf{i}_{obs}^{(h)}$, denoted by $\x_{obs}^{(h)}$;
\item the variables other than $\X^{(h)}$ and the observations whose indexes belong to $\mathbf{i}_{miss}^{(h)}$, denoted by $\x_{miss}^{(h)}$.
\end{itemize}

To begin, the original training data set $\sample$ is used to fill the blank spaces in a roughly way, for example, with the median of the observed values in the variable. Then, for each variable $\X^{(h)}$ a random forest is trained with target $\y_{obs}^{(h)}$ and predictors $\x_{obs}^{(h)}$, then the missing values $\y_{miss}^{(h)}$ are imputed with the prediction of $\x_{miss}^{(h)}$ using the random forest.

\subsubsection{Without imputation}
\paragraph{Missing Incorporated in Attributes.}

The Missing Incorporated in Attributes (MIA) consists in keeping all the missing values together when a split is performed. That is, missing values are assigned together to the child node that maximizes the CART criterion (or any other considered criterion). Thus, the splits with this approach assign the values according to one of the following rules:

\begin{itemize}
\item $\{\X^{(h)}<z\text{ and }\M^{(h)}=1\}$ versus $\{\X^{(h)}\geq z\}$,
\item $\{\X^{(h)}<z\}$ versus $\{\X^{(h)}\geq z\text{ and }\M^{(h)}=1\}$,
\item $\{\M^{(h)}=0\}$ versus $\{\M^{(h)}=1\}$.
\end{itemize}

%

\subsection{Description of the parameters}
\label{sec:simulations}

The regression function considered in this study is the so-called ``friedman1'' \citep{friedman1991multivariate}, which has been used in previous simulation studies \citep{friedman1991multivariate,breiman1996bagging,rieger2010random,josse2019consistency,friedberg2020local}, given by \[m(\x)=10\sin\parent{\pi\x^{(1)}\x^{(2)}}+20\parent{\x^{(3)}-0.5}^2+10\x^{(4)}+5\x^{(5)}\]
Our simulation study is based on the previous work of \citet{rieger2010random}, with the following characteristics:

\begin{itemize}
\item We simulate $\X$ uniformly distributed on $[0,1]^5$ and introduce missing values in $\X^{(1)},\X^{(3)}$ and $\X^{(4)}$, considering 7 different missing-data mechanisms.
\item For each missing-data mechanism we create 100 training data sets, each one with 200 observations.
\item In $\X^{(1)}$ 20\% of the data is missing, in $\X^{(3)}$ the amount is 10\%, and in $\X^{(4)}$ there is 20\% again.
\item We also create a testing data set with 2000 observations without missing values.

This amount of data is to have an appropriate approximation to the mean squared error (MSE) \[\E_{\X|\sample}\corchet{m_{M,n}(\X)-m(\X)}^2\] and the bias \[\E_{\X|\sample}\corchet{m_{M,n}(\X)-m(\X)}.\]
Note that these expressions are conditioned on the training sample $\sample$ and thus they are random variables which take a different value for each one of the 100 training data sets.

\item A random forest is built for each training data set and each missing-data mechanism as well as for the data sets without missing values (which are used as benchmark).
\item We use $M=100$ trees, which has been seen by simulation to be sufficient to stabilize the error in the case of the complete data sets.

For the rest of parameters we use the default values in the regression mode of the \texttt{R} package \texttt{randomForests}.

\item The parameter $\mathtt{mtry}$ is set to $\lfloor p/3\rfloor$.
\item We have sampled without replacement, so $a_n$ is set to $\lceil 0.632 n\rceil$.
\item And $\texttt{nodesize}$ is set to 5.
\end{itemize}

These parameters are the same for the median-imputation and MIA approaches. For Breiman's approach, Ishioka's approach and missForest we initialize the algorithms with the median-imputation and consider the same parameters as before for the construction of the regression trees. The number of iterations is set to 10 and random forests are built with 100 trees in each iteration, corresponding to the default values of the package \texttt{missForest} in \texttt{R}. We now describe the missing-data mechanisms, which are based on those presented by \citet{rieger2010random}.

\subsection{Missing-data mechanisms}
\subsubsection*{Missing at Random}

Methods for generating missing values at random are more complicated. The choice of the locations that are replaced by missing values in the ``missing'' variable now depends on the value of a second variable, called the ``determining'' variable. Therefore, the values of the ``determining'' variable now have influence on whether a value in the ``missing'' variable is missing or not. For $\X^{(1)}$ the ``determining'' variable is $\X^{(2)}$, while  $\X^{(5)}$ is used as the ``determining'' variable for $\X^{(3)}$ and $\X^{(4)}$.

\paragraph{Missing Completely at Random (MCAR).}
We select as many locations as desired sampled out of the $n$ observations and replace them by \texttt{NA}.

\paragraph{Creation of ranks (MAR1).} The probability for a missing value in a certain location in the ``missing'' variable is computed by dividing the rank of the location in the ``determining'' variable by $n(n+1)/2$. The locations for \texttt{NA} in the ``missing'' variable are then sampled with the resulting probability vector.

\paragraph{Creation of two groups (MAR2).} We divide the data set in two groups defined by the ``determining'' variable. A value belongs to the first group if the value in the ``determining'' variable is greater than or equal to the median of the ``determining'' variable, otherwise it belongs to the second group. An observation has a missing value with probability of 0.9 for the first group (0.1 for the second group) divided by the number of members in therespective group. The locations for \texttt{NA} in the ``missing'' variable are then sampled with the resulting probability vector.

\paragraph{Dexter truncation (MAR3).} The observations with the biggest values in the ``determining'' variable have the ``missing'' variable replaced by \texttt{NA} until the desired fraction of \texttt{NA} has been achieved.

\paragraph{Symmetric truncation (MAR4).} This method is similar to the previous one but we replace by \texttt{NA} the values in the ``missing'' variable in the observations with the biggest and the smallest values in the ``determining'' variable.

\paragraph{Missing depending on $Y$ (DEPY).} The missing values depend on the value of the response, the probability is 0.1 for observations where $Y\geq 13$, otherwise it is 0.4. The locations for \texttt{NA} in the ``missing'' variable are then sampled with the resulting probability vector.

\subsubsection*{Missing Not at Random}

\paragraph{Logit modelling (LOG).} In this method the probability for \texttt{NA} no longer depends on a single ``determining'' variable but in all the other variables. It is modeled as
\[\text{logit}\parent{\P\corchet{\M^{(h)}=1}}=-0.5+\sum_{\begin{tinymatrix} k=1 \\ k\neq h \end{tinymatrix}}^5 \X^{(k)}\]
Therefore, the probability of missingness depends on variables with observed values and variables with missing values.

\subsubsection*{Complete Observations}

Additionally, we consider the data sets with no missing values as a benchmark, which is denoted as ``COMP''.

\section{Results}
\label{sec:results}

\subsection{Change of missing-data mechanism}

\begin{figure}[H]
\begin{center}
\includegraphics[width=\textwidth]{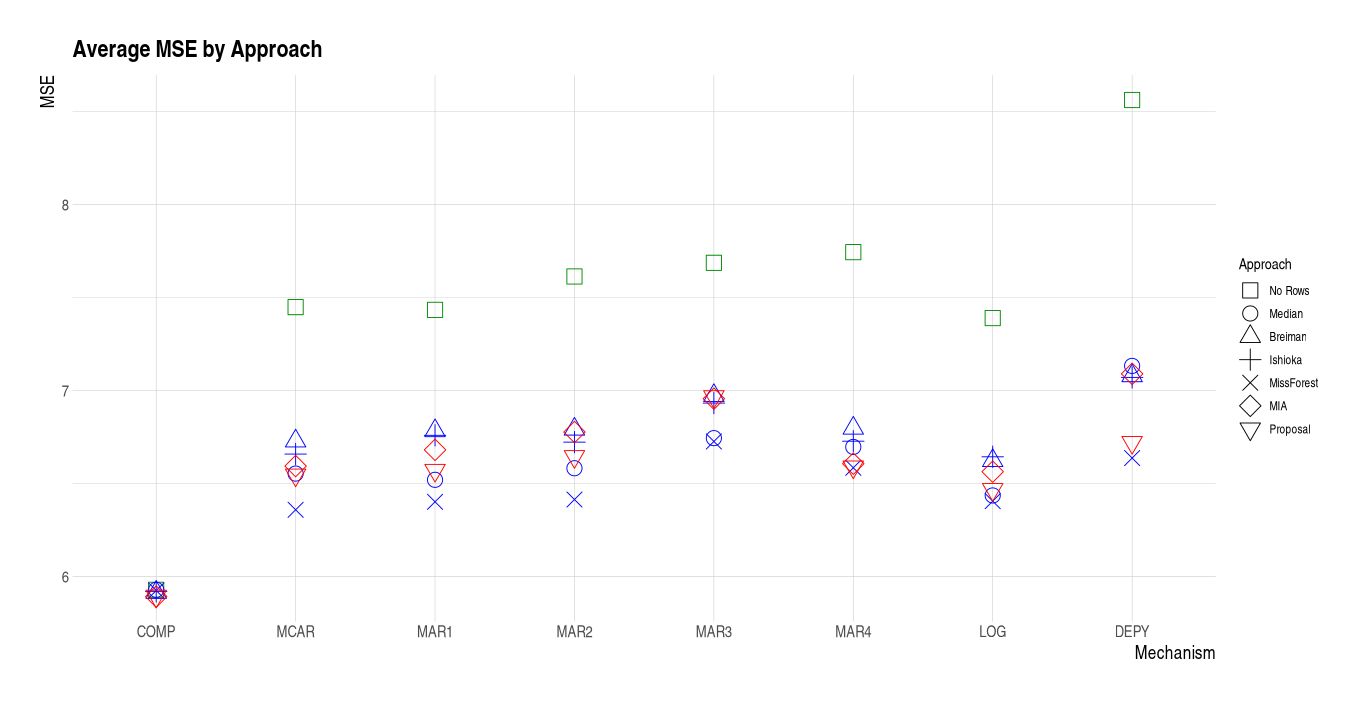}
\caption{Average MSE  of the testing data set for each approach and each missing-data mechanism.}
\label{fig:MSE_for_all_approaches}
\end{center}
\end{figure}

\begin{figure}[h!]
\begin{center}
\includegraphics[width=\textwidth]{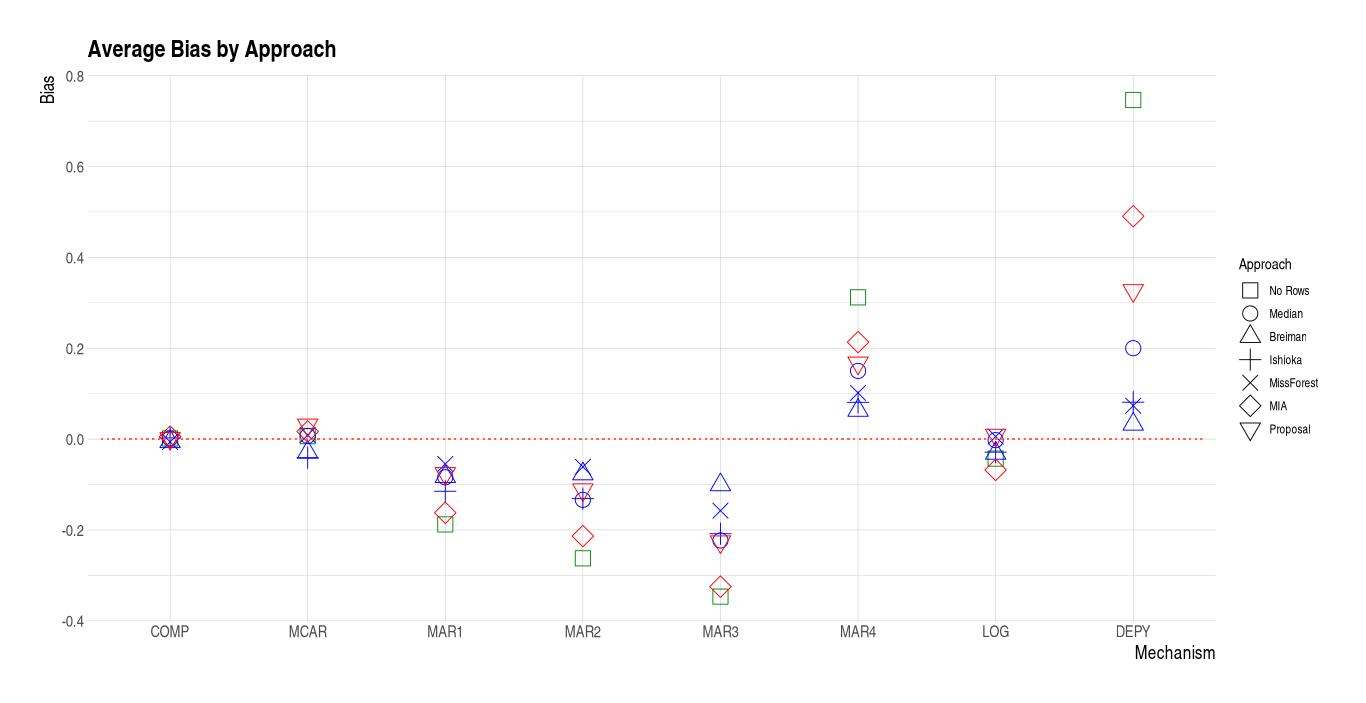}
\caption{Average bias  of the testing data set for each approach and each missing-data mechanism.}
\label{fig:Bias_for_all_approaches}
\end{center}
\end{figure}

\Cref{fig:MSE_for_all_approaches,fig:Bias_for_all_approaches} present the average MSE and the average bias over all the training data sets for each of the approaches considered. In green there is listwise deletion, those approaches that implement some imputation in the data set (Breiman's approach, Ishioka's approach and missForest) are in blue and those approaches that handle missing values directly in the construction of the trees (MIA and the our proposal) are in red. Listwise deletion (denoted as ``NoRows'') generates the largest errors and the estimates with more bias for all the mechanisms. Hence, it could be taken as a bound of the minimum expected performance for a method that attempts to estimate the regression function with missing values. We observe that missForest consistently generates estimators with the lowest errors regardless of the missing-data mechanism. Furthermore, we observe that our proposed method outperforms MIA, Breiman's approach and Ishioka's approach and can achieve similar errors as missForest. It is worthy to observe that even a simple approach as imputing with the median can outperforms most of the methods considered or with similar behavior in several missing-data mechanisms (see MAR1, MAR2, MAR3, LOG). In terms of bias, we observe that the algorithms that impute the missing values before the construction of the random forest tend to generate less biased results, while MIA tends to generate the second more biased estimators (with listwise deletion generating the most biased estimations). For the MCAR case, all the methods considered tend to be unbiased and with similar MSE, except for missForest and listwise deletion with the lowest and highest MSE, respectively.\footnote{Codes to reproduce our results can be found in: \url{https://github.com/IrvingGomez/RandomForestsSimulations}}

\subsection{Increasing the rate of missingness}
\label{sec:Rate_missingness}

\begin{figure}[H]
\begin{minipage}{0.45\textwidth}
\begin{center}
\includegraphics[width=\textwidth]{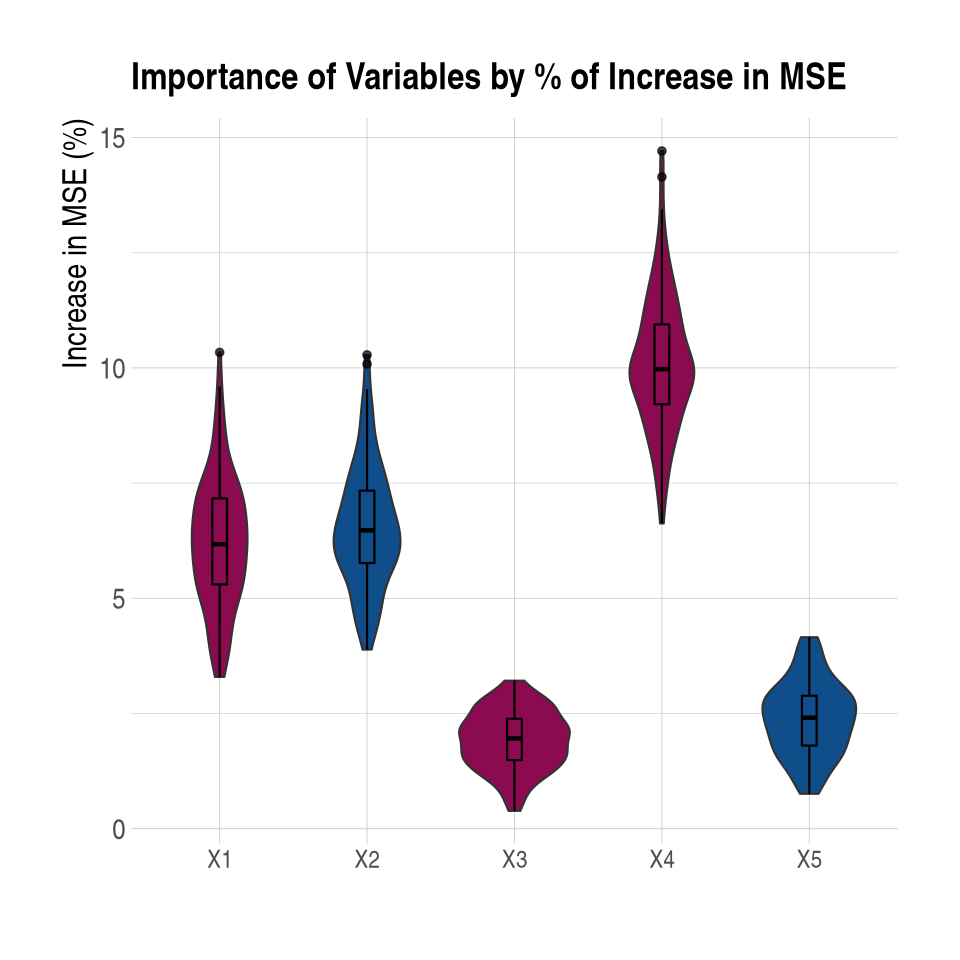}
\caption{Importance Variable accordingly to percentage increase in MSE.}
\label{fig:Importance1}
\end{center}
\end{minipage}\hspace{1cm}
\begin{minipage}{0.45\textwidth}
\begin{center}
\includegraphics[width=\textwidth]{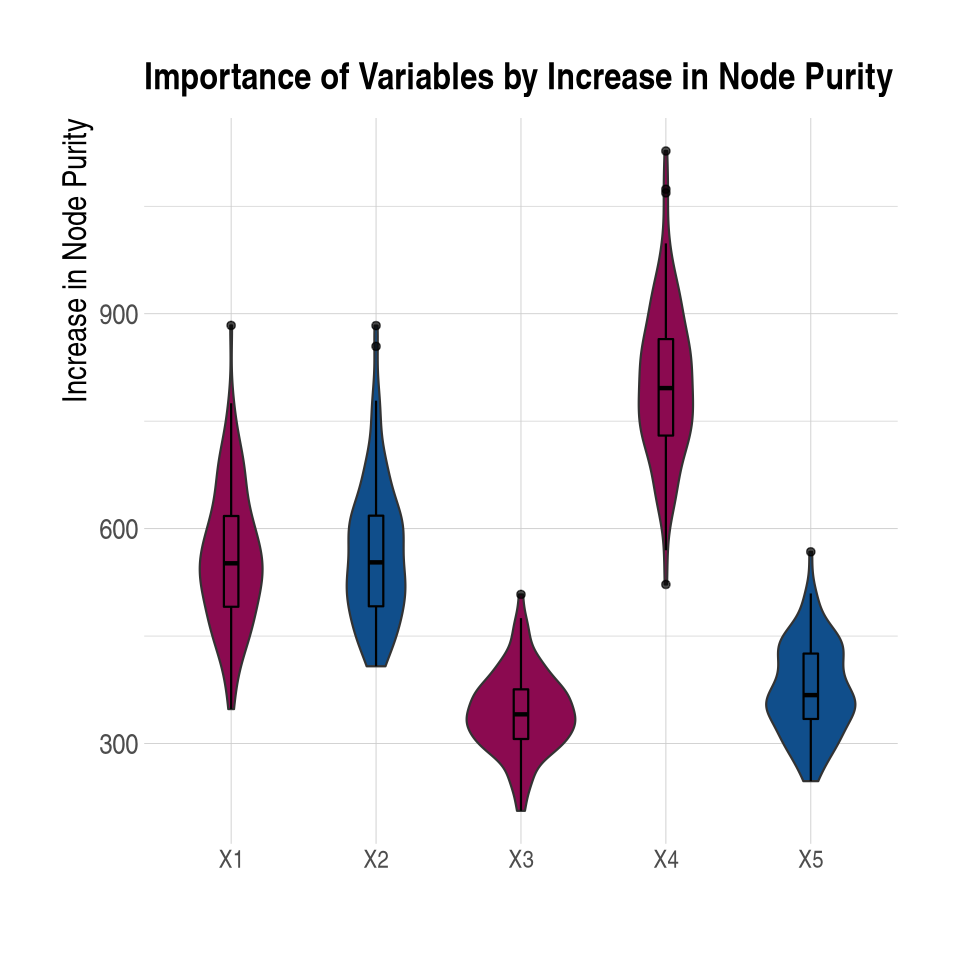}
\caption{Importance Variable accordingly to increase in node purity.}
\label{fig:Importance2}
\end{center}
\end{minipage}
\end{figure}

Using the 100 training data sets with no missing values, we calculate the importance of the variables with the \texttt{R} package \texttt{randomForests}, by percentage of increase in mean squared error and by increase in node purity \citep{breiman2001random, Breiman2003}. \Cref{fig:Importance1,fig:Importance2} show the violin plots for these measurements of importance. We can observe a consistently order for the variables with missing values in both measures of importance, where $\X^{(4)}$ is considered more important than $\X^{(1)}$ and $\X^{(3)}$. Hence, we decided to change the fraction of missingness in $\X^{(4)}$ to vary between 5\%, 10\%, 20\%, 40\%, 60\%, 80\%, 90\% and 95\%, without changing the percentage of missingness in $\X^{(1)}$ and $\X^{(3)}$. In this part of the study we do not consider anymore listwise deletion.

\begin{figure}[h!]
\begin{center}
\includegraphics[width=\textwidth]{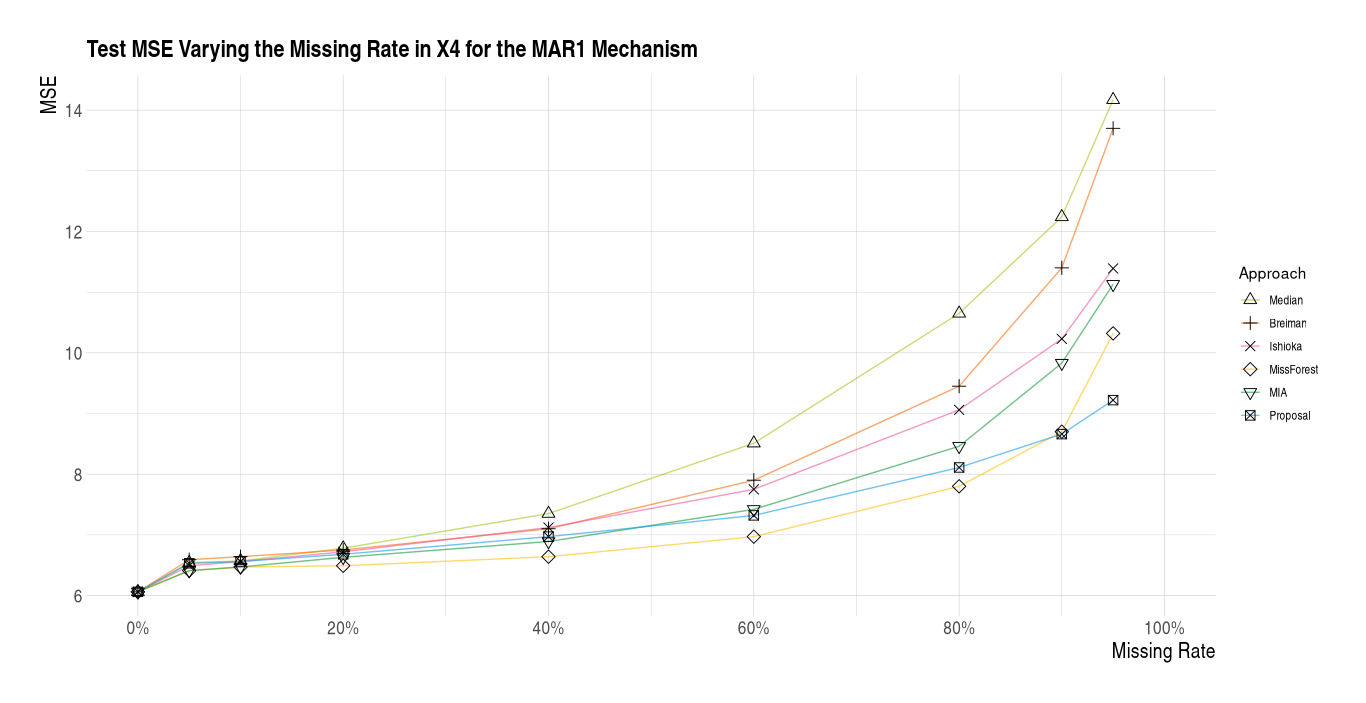}
\caption{Average MSE for the testing data set for each percentage of missingness, considering the MAR1 mechanism.}
\label{fig:MSE_changing_missingness}
\end{center}
\end{figure}

\Cref{fig:MSE_changing_missingness} presents the average MSE varying the percentage of missing data and considering the MAR1 mechanism\footnote{Analogous figures for the MSE and bias for all the missing-data mechanisms can be found in the supplementary material at \url{https://irvinggomez.com/publication/supplementary_random_forests_simulation/Supplementary_RandomForestsSimulation.pdf}}. Important differences between the performance of the methods become clear with the increasing percentage of missingness. Especially when this value is over 60\% there is an order on the performance of the methods. In this cases missForest, our proposal and MIA represent the methods with less MSE. Moreover, we can see the advantage of searching for the best assignation when the percentage becomes large with our proposal outperforming all the other algorithms. While all the methods present a similar MSE when the percentage of missing values is smaller than 40\%.

The differences between the missing-data mechanisms are also reflected when we increase the percentage of missing values. For example, the top three algorithms (missForest, our approach and MIA) present a MSE between 7.95 and 8.19 when the percentage of missingness is 90\% and the introduction of missing values is done completely at random (MCAR). On the other hand, for the same percentage of missing values and the same algorithms we observe a deterioration in terms of the MSE which varies between 9.05 and 12.86 when we consider the DEPY scenario. In this case we consistently observe that our approach and missForest outperform the other methods, regardless of the percentage of missing values, while there is no clear advantage for the rest of the algorithms over the others.

\begin{figure}[h!]
\begin{center}
\includegraphics[width=\textwidth]{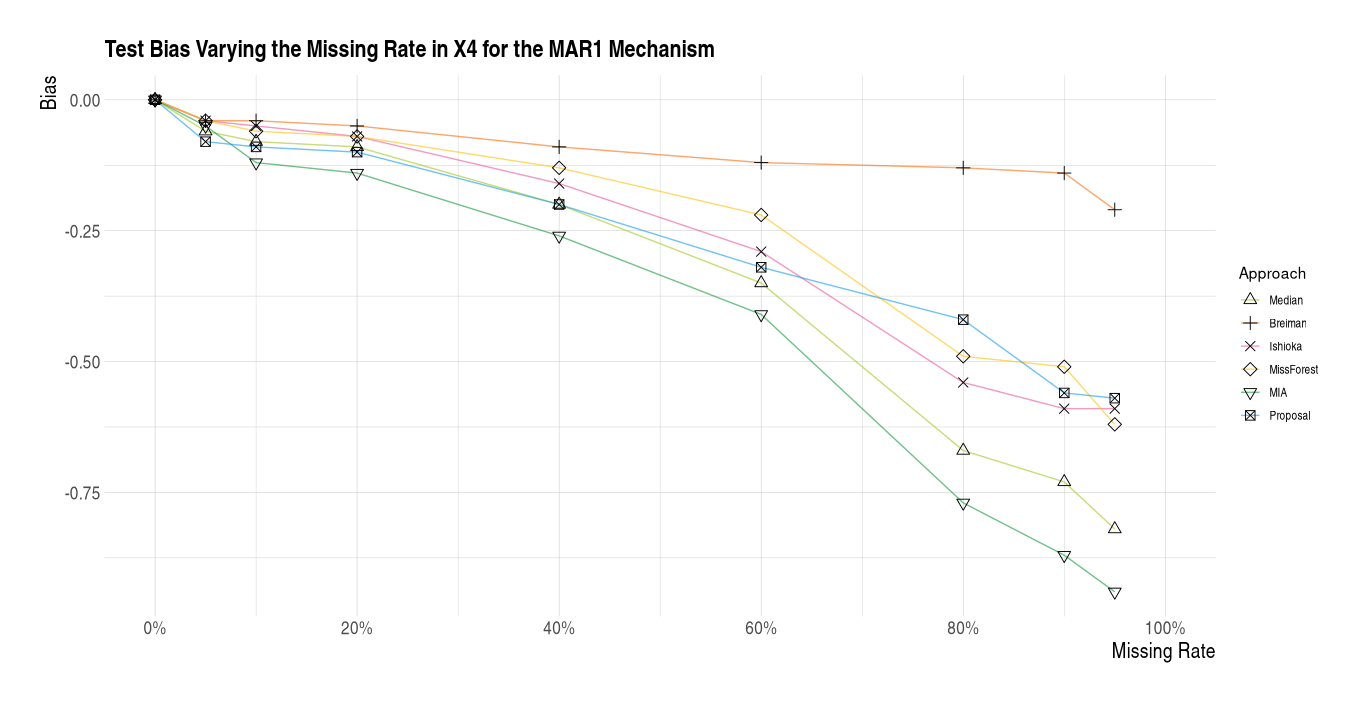}
\caption{Average bias for the testing data set for each percentage of missingness, considering the MAR1 mechanism.}
\label{fig:Bias_changing_missingness}
\end{center}
\end{figure}

\Cref{fig:Bias_changing_missingness} presents the average bias varying the percentage of missing data and considering the MAR1 mechanism. When we include the bias in the study, the differences between methods and missing-data mechanisms become more evident.  We observe that this missing-data mechanism introduced a bias in the methods. MIA and median-imputation create the most biased estimators, while Breiman's approach creates the less biased. Considering other missing-data mechanisms, we observe an unbiased behavior for the MCAR and the LOG cases, with some deterioration when the percentage of missing values affects up to 60\% of the observations.

\subsection{Prediction of new observations with missing entries}
\label{sec:prediction_with_missing}

A challenge in machine learning is to be able to compute a solution or a prediction to a certain problem when the given information is somewhat incomplete. In our specific context, our proposal do not have to rely on imputations of the missing entries to be able to calculate a prediction for a new data point. The fact that most of the existing techniques use imputation to construct the random forest tend to provide a prediction that is highly dependent on those imputations. In the following, we explain a way to perform the testing phase without having to impute the missing entries of the new data point.

\begin{figure}[H]
	\includegraphics[width=\textwidth]{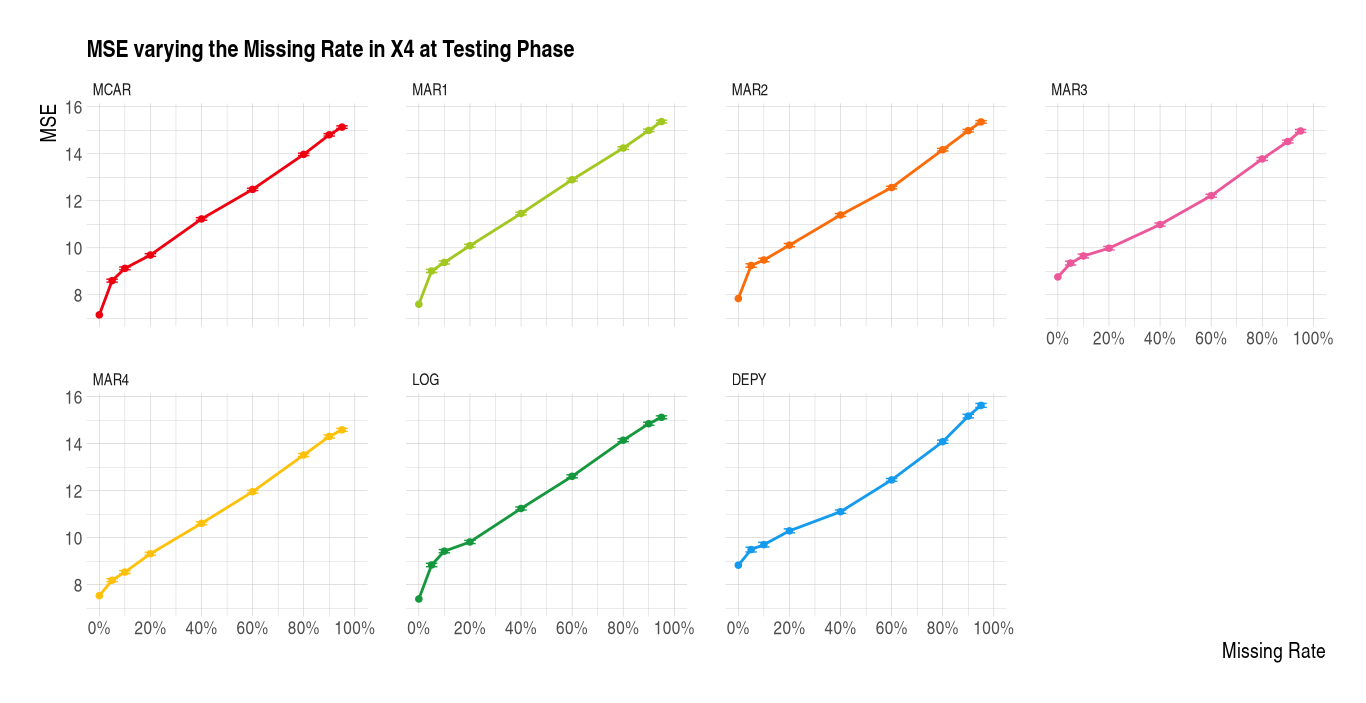}
	\caption{We compute the MSE over a testing set of 2000 data points, the proportion of missing values for \(\X^{(4)}\) varies from 0\% to 95\%.}
	\label{vary_testing}
\end{figure}

At this stage, we assume that the training phase is finished and that one have access to the random forest with the assignations associated to each tree.
The prediction is computed tree by tree and is averaged over the different trees in the random forest.
A tree prediction is performed looking for a (pseudo-random) final cell in the tree that is the most likely to contain the query point \(\X\). This is done in a recursive manner by going down in the tree following a procedure that we describe now.
Assume that we are in the cell \(A\) which has a cut $(h,z)$ that splits it into \(A_L\) and \(A_R\) and that the assignation vector associated to that cut is given by \(w\). If the direction $h$ is observed then, as usual, $\X$ is assigned to the left if $\X^{(h)}<z$, otherwise it is assigned to the right. Let us now assume that the direction $h$ is missing, in this case we need to look at the vector $w$ to assign the query point. Once again, let $N$ be the number of points with a missing value in the direction $h$ let $N_L$ be (resp. $N_R=N-N_L$) the number of such points assigned to the left (right) node. At this step, if \(N=0\) (when no missing values where observed in the direction \(h\) in the cell during the training phase) we stop the descent and predict the value of \(\Y\) by its mean value in the cell \(A\). Otherwise, \(N\neq 0\) and we can compute $p_L=N_L/N$ (resp. $p_R=N_R/N$) the probability for a missing observation to belong to the left (resp. right) node, given that it belongs to cell \(A\).
The next cell is, then, stochastically selected to the left or to the right with respective probabilities $p_L$ and $p_R$.
Our algorithm keeps track of the assignations \(w\) at each step so that computing the probabilities $p_L$ and $p_R$ is direct. We can think of the same procedure for the other random forest techniques that do perform a certain kind of assignation. For example, MIA algorithm is the most suited for a comparable treatment for incomplete data points in the testing phase. Nevertheless, the ``assignation information'' is not accessible in the existing codes for MIA. For the case of the algorithms that do perform imputation of the missing values, it is not clear how an imputation has to be done for the testing phase. As discussed in \citet{gomezmendez2020consistency}, these algorithms  are really dependent on the way the imputation is performed and then no theoretical guaranties are known for their consistency.
In the following, we give some simulations of this testing phase letting the proportion of missing values to vary. \\
In the training phase:
\begin{itemize}
	\item For each data set and each of the 7 missing-data mechanisms  (MCAR, MAR1, MAR2, MAR3, MAR4, LOG y DEPY) we introduce missing values in the variables \(\X^{(1)}\), \(\X^{(3)}\) and \(\X^{(4)}\). The proportion of missing values are 20\% for \(\X^{(1)}\), 10\% for \(\X^{(3)}\) and 60\% for the variable \(\X^{(4)}\).
	\item The random forests are constructed following the lines of Section \ref{sec:simulations}
\end{itemize}
For the testing phase:
\begin{itemize}
	\item For each missing-data mechanism, we let the percentage of missing values for \(\X^{(4)}\) to vary from 0\% to 95\% and let the other two proportions of missing values for \(\X^{(1)}\) and \(\X^{(3)}\) unchanged (at 20\% and 10\%, resp.)
	\item The graphics in Figure \ref{vary_testing} represent the MSE on those 2000 observation where 100 random forests are computed for each missing-data mechanism.
\end{itemize}

\section{Discussion and conclusions}
\label{sec:conclusions}

We developed a simulation study comparing 7 distinct approaches based on random forests to perform regression with missing entries. The algorithm proposed in this work handle missing values directly in the split criterion considered in the construction of the trees, assigning the each missing value to the node that maximizes the split criterion. Three of these algorithms use the random forests to impute the missing values and create completed data sets. Two of these group of algorithms rely on the computation of the proximity matrix and improve the imputations iteratively. The third algorithm, corresponding to missForest performs imputation through the implementation of random forests where the imputation of the missing values is treated as a regression problem by itself. We also considered MIA which, similar to our proposed algorithm, handle missing values directly in the construction of the trees. Finally, as simple benchmarks we considered the median-imputation and listwise deletion.

With no surprise, listwise deletion was the approach with the worst performance, this method should be avoided unless the percentage of observations with missing values is so low that they can be deleted without a severe harmful. For the rest of the algorithms, we observed differences between distinct techniques. These differences become more evident when the percentage of missing values increases, especially when it is over 60\%, while these differences are diluted for small values of this percentage (less than 40\%). Moreover, the behavior of the methods appear to be dependent on the missing-data mechanism. Intensive computer algorithms seem to perform particularly well for large percentage of missing values. In particular, we see that simple techniques as median-imputation have similar performance to more complicated algorithms when the percentage of missing values is under 40\% which makes its use sufficient in practice for few missing value contexts.

We studied the complexity of the proposal in \Cref{sec:algo_complexity} and show in \Cref{sec:simplifications} how a bisection technique can be performed to find the best assignation of the missing values. With these simplifications, we proved that the complexity of the algorithm is comparable to the MIA algorithm complexity up to a logarithmic factor, making possible to apply the algorithm for real applications.

Finally, in \Cref{sec:prediction_with_missing} we explain a process that allows the prediction of a new observation with missing entries. These procedure uses the assignation of the missing values during the training phase to estimate the probabilities of belonging to each cell. Hence, the new observation can be assigned stochastically to the cells of each tree using these probabilities. In the end, the prediction of the random forests is simply the average of the prediction of each tree. We observed through the simulation study that these technique seems to be robust to the missing-data mechanism and that it could be applied even for data sets with several missing values. Moreover, it can be applied with any other technique that assigns the missing observations, as long as the information of these assigantions is kept.


\printbibliography

\end{document}